\documentclass[12pt]{amsart}

\newcommand{\note}[2][\null]{%
 \marginpar{ 
   \renewcommand{\baselinestretch}{.5} 
   \vspace{-1em}\hrule\vspace{3pt}%
   \footnotesize\raggedright\textsf{#2\ifx#1\null\else\\\hfill--- 
     {\em #1}\fi}\vspace{1.5em} 
 }%
} 

\usepackage{amsmath}
\usepackage{mathrsfs}
\usepackage{amsthm}
\usepackage{graphicx}
\usepackage{amsfonts}
\usepackage{wasysym}
\usepackage{textcomp}
\newcommand{\re}{\mathrm{Re}}

\begin{document}
\title{On the Minkowski-Funk Transform}
\author{Susanna Dann}

\address{Mathematics Department\\
Louisiana State University\\
Baton Rouge, Louisiana}
\email{sdann@math.lsu.edu}

\begin{abstract}
The subject of this paper is the history of the Minkowski-Funk Transform. After introducing the Minkowski-Funk Transform as well as its dual transform and a generalization of both, we will present an inversion formula of the Minkowski-Funk Transform. Then we will discuss the history of this problem: related work by Minkowski and Funk and the connection between their work.
\end{abstract}
\maketitle

\section{Notation and Preliminaries}
Let $\mathbb{R}^n$ denote the real $n-$dimensional Euclidean space, $\{e_1, e_2, ..., e_n \}$ the standard basis for $\mathbb{R}^n$, $S^{n-1} \subset\mathbb{R}^n $ the unit sphere, and $\sigma_{n-1}=\frac{2 \pi^{n/2}}{\Gamma(n/2)}$ the surface area of  $S^{n-1}$. Let $\Xi$ be the set of $(n-2)-$dimensional totally geodesic submanifolds or $(n-2)-$geodesics $\xi \subset S^{n-1}$. Each  $(n-2)-$geodesic is a section of  $S^{n-1}$ with a hyperplane in  $\mathbb{R}^n$ through the origin. If $n=3$, then $(n-2)-$geodesics correspond to the great circles. For $\theta \in S^{n-1}$, let $\theta^{\perp}$ denote the hyperplane through the origin having $\theta$ as its normal vector, i.e. $\theta^{\perp}= \{x \in \mathbb{R}^n : <x,\theta> = 0\}$, where $ <.,.> $ denotes the standard inner-product on $\mathbb{R}^n$. Let $\xi_o = S^{n-1} \bigcap e_n^{\perp} $. For $\theta \in [0, \frac{\pi}{2}]$, let $g_{k,n}(\theta)$ denote the rotation in the $e_k e_n-$plane with the matrix 
$ \left( \begin{array}{cc} \sin \theta & \cos \theta \\ - \cos \theta & \sin \theta \end{array} \right)$, which is a clockwise rotation by the angle $\frac{\pi}{2}-\theta$. Further, let $x_{\theta}$ denote the image of $e_n$ under this rotation: $x_{\theta} = g_{k,n}(\theta) e_n = e_k \cos \theta + e_n \sin \theta$. Note that $d(x_\theta, \xi_o)=\theta$, where $d:S^{n-1} \times S^{n-1} \rightarrow [0, \pi] $ is the geodesic-distance function on the unit sphere.
\\
Let $G=SO(n)$, $K=SO(n-1)$. $G$ acts transitively on $S^{n-1}$ and $\Xi$, while $K$ is the stabiliser of $e_n$ and $\xi_o$, hence $S^{n-1} \cong G/K$ and $\Xi \cong G/K$. For $\theta \in S^{n-1}, \xi \in \Xi$, let $r_{\theta}, r_{\xi}$ denote rotations such that $r_{\theta} e_n = \theta, r_{\xi} \xi_o = \xi$. Since $G$ is a compact Hausdorff topological group, we have a unitmodular Haar measure on $G$. Let $\varphi$ be a measurable function on $\Xi$. We define a $G-$invariant measure $d\xi$ on $\Xi$ by $$ \int_{\Xi} \! \varphi(\xi) d\xi =  \int_G \! \varphi(g\xi_o) dg,$$ where $\int_G \!  dg = 1$.
For $\theta \in S^{n-1}$ and $t \in (-1, 1)$ let $S_{\theta, t} = \{\sigma \in S^{n-1}: \sigma \cdotp \theta = t \}$, an $(n-2)-$dimensional sphere of radius $\sqrt{1 - t^2}$, and $$ (M^t f)(\theta) = c_n  \int_{S_{\theta, t}} \! f(\sigma) d_{\theta}\sigma,$$ where $d_{\theta}\sigma$ denotes the Lebesgue measure on $S_{\theta, t}$ and $c_n$ is the normalizing factor: $c_n = \frac{(1 - t^2)^{\frac{2-n}{2}}}{\sigma_{n-2}}$. $M^t$ is an average of a function $f$ over $S_{\theta, t}$.\\
\textit{Remark:} $M^t$ commutes with rotations.\\
\\
The Riemann-Liouville fractional integrals and derivatives on $(0, \infty)$ will be defined as in (\cite{rub06}, p. 16 -18) by 
$$ (I_{0+}^{\alpha} \psi)(t) = \frac{1}{\Gamma(\alpha)}  \int_0^t \! \psi(\tau) (t-\tau)^{\alpha-1} d\tau, \; \; \; \; \; D_{0+}^{\alpha} = (I_{0+}^{\alpha})^{-1}.$$
\\
The \textsl{\textbf{ Minkowski-Funk transform}} of a function $f$ on $S^{n-1}$ is a function $Mf$ on $\Xi$ given by 
\begin{equation}\label{eqm} (Mf)(\xi)= \int_{\xi} \! f(x) \, dx. \end{equation}
For this definition to make sense we need to fix the measure on the $(n-2)-$geodesics in a coherent way. Define
$$ \int_{\xi} \! f(x) \, dx = \int_{\xi_o} \! f(r_{\xi} \theta) \, d\theta,$$   where $d \theta$ denotes the Lebesgue measure on $\xi_o=S^{n-2}$ and $\int_{\xi_o} \! d\theta = \sigma_{n-2}$. Furthermore, using repeated slice integration we get
$$\int_{\xi_o} \! f(r_{\xi} \theta) \, d\theta = \int_0^{\pi} \int_0^{\pi} \dots \int_0^{\pi} \int_0^{2 \pi} \! f( r_{\xi} \theta(\varphi_1, \dots \varphi_{n-2})) $$ 

$$\sin^{n-3}\varphi_1 \sin^{n-4}\varphi_2 \dots \sin\varphi_{n-3} d\varphi_{n-2} d\varphi_{n-3} \dots d\varphi_2  d\varphi_1, $$ where $\varphi_1, \dots , \varphi_{n-2}$ are the spherical coordinates on $\xi_o$.\\
\\
The \textsl{\textbf{dual Minkowski-Funk transform}} of a function $\varphi$ on $\Xi$ is a function $M^{*}\varphi$ on $S^{n-1}$ given by $$ 
(M^{*}\varphi)(x)= \int_{x \in \xi} \! \varphi(\xi) \, d\xi.$$
We can generalize the transforms $f \rightarrow Mf$, $\varphi \rightarrow M^{*}\varphi $ as follows. Let $\theta \in [0, \frac{\pi}{2}] $, 
$$ (M_{\theta}f)(\xi)= \int_{d(x, \xi)=\theta} \! f(x) \, dm(x) = \int_K \! f(r_{\xi} \rho x_{\theta}) d\rho, $$
\begin{equation}\label{eqgm} (M_{\theta}^{*}\varphi)(x)= \int_{d(x, \xi)=\theta} \! \varphi(\xi) \, d\mu(\xi) = \int_K \! \varphi(r_x \rho g_{k,n}^{-1}(\theta) \xi_o) d\rho,  \end{equation}
where $dm, d\mu$ denote the normalized measures. For $\theta=0$, we get: $M_0 f = \sigma_{n-2}^{-1} Mf$, and $M_0^{*}\varphi = M^{*}\varphi$. Compare the above definitions of the Minkowski-Funk and related transforms with \cite{helg99}, \cite{rub02}.

\section{Inversion of the Minkowski-Funk transform} 
We have the following fact:
\textit{For $f \in L_{even}^1(S^{n-1})$ }
\begin{equation} (M_{\theta}^*Mf)(x) = 2 \pi^{\frac{n-2}{2}} \cos^{3-n} \theta (I_{0+}^{\frac{n-2}{2}} \widetilde{f_x}) (\cos^2 \theta),  \end{equation}
$$ \widetilde{f_x}(\tau) = \frac{1}{\sqrt{\tau}} (M^{\sqrt{\tau}} f)(x). $$
\textit{Proof:} Let $z \in S^{n-1}$, then
\begin{equation}\label{eq1} \int_K \! F(\rho z) d\rho = (M^{z_n}F)(e_n). \end{equation}
Set $z=g^{-1}\eta, \eta \in \xi_o, g=g_{n-2, n}(\theta), F(\cdotp)=f(r_x \cdotp)=:f_{r_x }(\cdotp)$, then $$lhs(\ref{eq1}) = \int_K \! f(r_x \rho g^{-1} \eta) d\rho,$$ and $z_n = z \cdotp e_n = g^{-1}\eta \cdotp e_n = \eta \cdotp g e_n = \eta \cdotp x_{\theta}$, so
$$rhs(\ref{eq1}) = (M^{\eta \cdotp x_{\theta}}f_{r_x})(e_n) = (M^{\eta \cdotp x_{\theta}}f)(r_x e_n) = (M^{\eta \cdotp x_{\theta}}f)(x), $$ where the first equality holds since $M^t$ commutes with rotations. So we get
\begin{equation}\label{eq2} \int_K \! f(r_x \rho g^{-1} \eta) d\rho = (M^{\eta \cdotp x_{\theta}}f)(x). \end{equation}
Now we integrate (\ref{eq2}) in $\eta \in S^{n-2} = \xi_o$ and interchange the order of integration. We get by (\ref{eqm}), (\ref{eqgm})
$$lhs(\ref{eq2}) = \int_K d\rho \int_{\xi_o} \! f(r_x \rho g^{-1} \eta) d\eta =  \int_K d\rho \int_{r_x \rho g^{-1} \xi_o} \! f(\eta) d\eta $$
$$ \!\!\! \!\!\! \!\!\! \!\! = \int_K \! (Mf)(r_x \rho g^{-1} \xi_o) d\rho = (M_{\theta}^*Mf)(x). $$
Since $x_{\theta} = e_{n-2} \cos \theta + e_n \sin \theta, \eta \cdotp x_{\theta} = (\eta \cdotp e_{n-2}) \cos \theta$, 
$$ rhs(\ref{eq2}) =  \int_{S^{n-2}} \! (M^{\eta \cdotp x_{\theta}}f)(x) d\eta = \int_{S^{n-2}} \! (M^{(\eta \cdotp e_{n-2}) \cos \theta}f)(x) d\eta.$$
Now we apply Catalan's Formula (\cite{rub09}, p. 5) to get
$$ rhs(\ref{eq2}) = \sigma_{n-3} \int_{-1}^1 \! (M^{\tau \cos \theta}f)(x) (1 - \tau^2)^{\frac{n-4}{2}} d\tau. $$
$\tau \mapsto  (M^{\tau \cos \theta}f)(x)$ is even, since $f$ is even (\cite{rub09}, p. 72), so
$$ rhs(\ref{eq2}) = 2 \sigma_{n-3} \int_{0}^1 \! (M^{\tau \cos \theta}f)(x) (1 - \tau^2)^{\frac{n-4}{2}} d\tau. $$
Changing the varialbes: $y = \tau \cos \theta$, we get
$$ rhs(\ref{eq2}) = 2 \sigma_{n-3} \int_{0}^{\cos \theta} \! (M^y f)(x) (\cos^2 \theta - y^2)^{\frac{n-4}{2}} \cos^{3-n} \theta dy. $$
Changing the varialbes again: $\tau = y^2$, we get
$$ rhs(\ref{eq2}) = 2 \sigma_{n-3} \int_{0}^{\cos^2 \theta} \! (M^{\sqrt{\tau}} f)(x) (\cos^2 \theta - \tau)^{\frac{n-4}{2}} \cos^{3-n} \theta \frac{d\tau}{2\sqrt{\tau}} $$
$$  = \sigma_{n-3} \cos^{3-n} \theta \int_{0}^{\cos^2 \theta} \!  \left\{ \frac{1}{\sqrt{\tau}} (M^{\sqrt{\tau}} f)(x) \right\} (\cos^2 \theta - \tau)^{\frac{n-2}{2}-1}  d\tau . $$
Setting $\widetilde{f_x}(\tau) = \frac{1}{\sqrt{\tau}} (M^{\sqrt{\tau}} f)(x)$, we obtain
$$ rhs(\ref{eq2}) = \sigma_{n-3} \cos^{3-n} \theta \; \Gamma(\frac{n-2}{2}) (I_{0+}^{\frac{n-2}{2}} \widetilde{f_x}) (\cos^2 \theta)$$
$$ = 2 \pi^{\frac{n-2}{2}} \cos^{3-n} \theta \; (I_{0+}^{\frac{n-2}{2}} \widetilde{f_x}) (\cos^2 \theta). $$
This gives the result. \qed \\
\\
Thus we have: $(I_{0+}^{\frac{n-2}{2}} \widetilde{f_x}) (\cos^2 \theta) = (2 \pi^{\frac{n-2}{2}})^{-1} \cos^{n-3} \theta (M_{\theta}^*Mf)(x).$ By changing the exterior variable on the left-hand-side, namely $t = \cos^2 \theta$, it becomes:
$$(I_{0+}^{\frac{n-2}{2}} \widetilde{f_x}) (t) = (2 \pi^{\frac{n-2}{2}})^{-1} t^{\frac{n-3}{2}} (M_{\cos^{-1}(\sqrt{t})}^*Mf)(x).$$
This is inverted by $D_{0+}^{\frac{n-2}{2}}$ (\cite{rub06}, Lemma 3.3 on page 18), thus
$$ \widetilde{f_x}(t) = (2 \pi^{\frac{n-2}{2}})^{-1} D_{0+}^{\frac{n-2}{2}} \left( (\cdotp)^{\frac{n-3}{2}} (M_{\cos^{-1}(\sqrt{\cdotp})}^*Mf)(x) \right)(t). $$
Finally, since $M^t f \rightarrow f$ as $t \rightarrow 1$ for $f \in L^1(S^{n-1})$ (\cite{rub09}, p. 69), $\lim_{t \rightarrow 1} \sqrt{t} \widetilde{f_x}(t) = f(x).$ We have followed the proof in \cite{rub02}.
 

\section{History of the Minkowski-Funk transform} 
\subsection{Funk's work}
Funk \footnote{Paul Georg Funk (14. April 1886 in Wien - 3. Juni 1969 in Wien) an Austian mathematician. He was a student of David Hilbert. He earned his PhD at Georg-August-Universit\"{a}t G\"{o}ttingen in 1911. The tilte of his thesis is: \"{U}ber Fl\"{a}chen mit lauter geschlossenen geod\"{a}tischen Linien".} has considered the following problem: Reconstruction of a function $f$ from its integrals over all great circles on the sphere $S^2$. It was part of his dissertation (1911) and in 1913 he published a longer paper \cite{funk13} (with the same title as his dissertation) that discussed this problem in chapter 2. In that paper (on page 284) for a function $\Phi$ on $S^2$  he defines a \textsl{\textbf{"circle-integral function of $\Phi$"}} \textsl{("die Kreisintegral-Funktion von $\Phi$") } on the set of great circles as an integral over the corresponding great circle, he denotes it by $\chi$. There he also derives a set of properties of $\chi$ (p. 284-285, 286) as well as two inversion algorithms. The first involves expansion in series of spherical harmonics (p. 285-286) and the second reduces the problem to the Abel's integral equation (p. 287-288). There he also indicates two geometric applications of the transform (p. 287). One of the applications beeing the derivation of the Minkowski's Theorem, \textit{that bodies of constant circumference are bodies of constant width}.

\subsection{Minkowski's work}
Minkowski \footnote{Hermann Minkowski (June 22, 1864 in Aleksotas (Russian Empire, now Lithuania) - January 12, 1909 in G\"{o}ttingen) was a German mathematician and physisist. He earned his PhD at the University of K\"{o}nigsberg in 1885. The title of his thesis is: "Untersuchungen \"{u}ber quadratische Formen. Bestimmung der Anzahl verschiedener Formen, welche ein gegebenes Genus enth\"{a}lt
".} proved the following fact related to the convex geometry in his paper from 1904 \cite{mink04}: \textsl{A body is of constant width iff it is a body of constant circumference}. To make this result presice we will need to make some definitions:\\
\\
Let $\Omega \subset \mathbb{R}^3$ be a compact convex body with smooth boundary containing the origin in its interior. A \textsl{\textbf{support plane}} is a plane that contains at least one point of $\Omega$, but does not intesect $\Omega$. Let $\omega \in S^2$, then $<x, \omega> = H(\omega)$ stands for the support plane perpendicular to $\omega$ at distance $H(\omega)$ from the origin. The distance between two support planes perpendicular to a given $\omega \in S^2$ is called the \textsl{\textbf{width $B(\omega)$}} of $\Omega$ in the direction $\omega$: $B(\omega) = H(\omega) + H(-\omega)$. A body with the same width in every direction is called a \textsl{\textbf{body of constant width}}. Next, consider the projection of $\Omega$ onto the plane through the origin perpendicular to $\omega$ and let $C(\omega)$ denote the boundary curve of this projection. The arc length of $C(\omega)$ is called the \textsl{\textbf{circumference $U(\omega)$}} of $\Omega$ in the direction $\omega$.  A body with the same circumference in every direction is called a \textsl{\textbf{body of constant circumference}}. Expanding the functions $B(\omega)$ and $U(\omega)$
in series of spherical harmonics Minkowski proved the above result.
\subsection{Connection between the Minkowski's result and \\ the  Minkowski-Funk transform}
The connection between the two is established through the following fact: \textit{Let $\Omega, \omega, B(\omega)$, and $U(\omega)$ be as defined above, then
$$
U(\omega)= \frac{1}{2}  \int_{S^2 \bigcap \omega^{\perp}} \! B(\sigma) ds(\sigma),
$$
where $ds$ is the arc-length element of the corresponding great circle.}
\\
\\
\textit{Proof:} 
The circumpherence of $\Omega$ in the direction of the $z-$axis is the length of a smooth closed curve in the $xy-$plane, which is the envelope of the following family of curves: $<p, \omega_{\varphi}>=H(\omega_{\varphi})$, where $0 \leq \varphi \leq 2\pi$. These are the projections of the vertical support planes onto the $xy-$plane, in other words, support planes with normal vectors parallel to the $xy-$plane. Any such normal vector $ \omega_{\varphi}$ has coordinates $<\cos \varphi, \sin \varphi>$, $0 \leq \varphi \leq 2\pi$. Let us write $h(\varphi)$ for the distance $H(\omega_{\varphi})$, which is the distance of the support plane, determined by the normal vector $\omega_{\varphi}$ and thus by the angle $\varphi$, to the origin. Further, let us write $F(x, y, \varphi)= <p, \omega_{\varphi}> - h(\varphi)$. $F(x, y, \varphi)=0$ is the family of curves (= support lines), which envelope the boundary curve $C_{e_3}$ of the projection of $\Omega$ onto the $xy-$plane. To write $C_{e_3}$ as a parametric curve $(x(\varphi), y(\varphi))$ one needs to solve the system of equations: 
$$ F(x, y, \varphi)=0,$$
$$ \frac{\partial F(x, y, \varphi)}{\partial {\varphi}} =0.$$
We obtain for $ 0 \leq \varphi \leq 2\pi $ :
$$C_{e_3}(\varphi) = (h(\varphi) \cos \varphi - \frac{\partial h}{\partial \varphi}(\varphi) \sin \varphi, h(\varphi) \sin \varphi + \frac{\partial h}{\partial \varphi}(\varphi) \cos \varphi).$$
The arc-length of a curve in a parametric form is computed by: $  \int \! \sqrt{(x^{'})^2 + (y^{'})^2}$. In our case:
$$ x^{'} = -(h + \frac{\partial^2 h}{\partial \varphi^2}) \sin \varphi, \;\;\;\;\;\;\;\;  y^{'} = (h + \frac{\partial^2 h}{\partial \varphi^2}) \cos \varphi. $$
So $(x^{'}(\varphi))^2 + (y^{'}(\varphi))^2 = (h + \frac{\partial^2 h}{\partial \varphi^2})^2$ and thus 
$$ U(e_3) =  \int_0^{2\pi} \! \mid h(\varphi) + \frac{\partial^2 h}{\partial \varphi^2} (\varphi) \mid d\varphi. $$
The integrand in the above integral is positive. To see this let us find another expression for it: $$ h + \frac{\partial^2 h}{\partial \varphi^2} = y^{'} \cos \varphi - x^{'} \sin \varphi = <\cos \varphi, \sin \varphi> <y^{'}, - x^{'}> . $$ So it is the dot product between the normals of the support line corresponding to $\varphi$ and of the tangent line to $C_{e_3}$ at the point $(x(\varphi), y(\varphi))$. Let us denote the angle between these two vectors by $\theta$. If $\theta = 0$, then the tangent line and the support line are the same. If we assume that $\theta \geq \frac{\pi}{2}$, we get that $0 \notin$ interior of $\Omega$, a contradiction. Thus $0 \leq \theta < \frac{\pi}{2}$. This implies $ h(\varphi) + \frac{\partial^2 h}{\partial \varphi^2} (\varphi) \geq 0$ for all $\varphi$. So 
$$U(e_3) =  \int_0^{2\pi} \! \left( h(\varphi) + \frac{\partial^2 h}{\partial \varphi^2} (\varphi) \right) d\varphi. $$
Since the boundary of $\Omega$ is smooth, $h(\varphi)$ is smooth as well, this gives $ \int_0^{2\pi} \! \frac{\partial^2 h}{\partial \varphi^2} d\varphi = \frac{\partial h}{\partial \varphi} \mid_0^{2\pi} = 0$ and 
$$ \int_0^{2\pi} \! h(\varphi) d\varphi = \int_0^{2\pi} \! H(\omega_{\varphi}) d\varphi = \int_0^{2\pi} \! H(-\omega_{\varphi}) d\varphi. $$
This gives the result. \qed \\
\\
This result can be found in the Minkowski's original paper (\cite{mink04}, p. 506-507) and in (\cite{helg07}, p. 129). So in Minkowski's work we see an integral over a great circle, where the function integrated and the resulting function both have a geometric meaning. Funk, assigned by his advisor to look at Minkowski's work, considered then such types of integrals for arbitrary functions on the sphere.  \\
\\
Remark to the above proof: There is also another way to verify that $  h(\varphi) + \frac{\partial^2 h}{\partial \varphi^2} (\varphi) \geq 0 $. Let $ \kappa $ denote the curvature of a curve, then for a plane curve $ \mathscr{C}$ given by a parametric equation $(x(t) , y(t))$ we have: $ \kappa = \frac{|x^{'}y^{''}-y^{'}x^{''}|}{(x^{'2}+y^{'2})^{\frac{3}{2}}}$. If one takes into account the direction in which the tangent vector rotates as one moves along the curve, then one considers the signed curvature. Namely, if the tangent vector rotates counterclockwise, then $\kappa > 0$, and if it rotates clockwise, then $ \kappa < 0$. \\
In our case the curvature of the boundary curve $ C_{e_3}(\varphi) $ is $ \kappa(\varphi) = \frac{1}{h(\varphi) + \frac{\partial^2 h}{\partial \varphi^2} (\varphi)}$ and the sign of the curvature is positive. This gives us exactly what we want. 

\subsection{Solution of the Minkowski's problem be means of the  Minkowski-Funk transform}
We will use the following fact (\cite{helg07}, p. 125): \textit{The kernel of the Minkowski-Funk transform consists of odd functions.}\\
Assuming that $B(\omega) = c$ for every $\omega \in S^2$, where $c$ is some constant, we get $U(\omega) = \pi c$. On the other hand, assume  $U(\omega) = c$ for every $\omega \in S^2$. Writing $c = \frac{1}{2\pi}  \int_{S^2 \bigcap \omega^{\perp}} \!  c \; ds(\sigma)$, we get $\int_{S^2 \bigcap \omega^{\perp}} \! \left( \frac{B(\sigma)}{2} - \frac{c}{2\pi} \right) ds(\sigma) = 0$. Thus, by the above fact $ \frac{B(\sigma)}{2} - \frac{c}{2\pi}$ is odd, which implies that $B(\sigma)=\frac{c}{\pi}$ for every $\sigma \in S^2$.
So we obtain $U \equiv$ constant iff  $B \equiv$ constant. \qed \\
\\
There is another way to obtain Minkowski's theorem from the Minkowski-Funk transform, namely using the theory of spherical harmonics. For this we recall some definitions and facts, compare \cite{rub09} for complete statements. \\
For $\alpha$ with $\re(\alpha) \geq 0$ and $\alpha \neq 1, 3, 5, ...,$ let $$\left( \widetilde{M}^{\alpha}f \right)(\theta) = \gamma_n(\alpha) \int_{S^{n-1}} \! f(\sigma) |\theta \cdotp \sigma|^{\alpha-1} d\sigma. $$
The Fourier-Laplace decomposition of $\widetilde{M}^{\alpha}f$ for $f \in C^{\infty}(S^{n-1}),$ has the form $$ \widetilde{M}^{\alpha}f = \Sigma_{m, \mu} h_m(\alpha) f_{m, \mu} Y_{m, \mu}.$$ 
If f is an even continuous function and  $\alpha + \beta = 2-n$, then $ \widetilde{M}^{\alpha} \widetilde{M}^{\beta} f = f$ and $\lim_{\alpha \rightarrow 0} \left( \widetilde{M}^{\alpha}f \right)(\theta) = (2 \pi^{\frac{n-2}{2}})^{-1} \left( Mf \right)(\theta).$ \\
\\
We will show: \textit{ $f \equiv$ constant iff $Mf \equiv$ constant.}\\
Assume $f \equiv constant$, then $f$ is homogeneuos of degree $0$, thus the Fourier-Laplace coefficients of $f$, $f_{m, \mu} = 0$ for $m > 0$. $d_n(0)=1$, so $$ \left( Mf \right)(\theta) = h_0(0) f_{0,0} Y_{0,0}(\theta) \equiv constant.$$
Conversely, assume $  \left( Mf \right)(\theta) = \varphi \equiv constant$, then $$f=\widetilde{M}^{2-n} \varphi = \Sigma_{m, \mu} h_m(2-n) \varphi_{m, \mu} Y_{m, \mu}.$$ Again, since $\varphi$ is constant, only $\varphi_{0,0} \neq 0$. So $f = h_0(0) \varphi_{0,0} Y_{0,0}(\theta) \equiv constant.$ \qed \\
\\
Since $ U(\omega) = \frac{1}{2} \left(MB\right) (\omega)$, the Minkowski's theorem follows.

\section*{Acknowledgements}
This paper is based on the project for the class \textit{Introduction to the Theory of Spherical Harmonics}.
The author thanks Professor Dr. Rubin for suggesting this topic for the class paper and for helpful discussions.

\end{document}